\documentclass[12pt]{article}
\usepackage{e-jc}

\title{The Tournament Scheduling Problem with Absences\footnotetext{\copyright\ 2016. This work is licensed under the Creative Commons
Attribution-NonCommercial-NoDerivatives 4.0 
License
http://creativecommons.org/licenses/by-nc-nd/4.0/
}
\footnotetext{Published as: U. Schauz, The tournament scheduling problem with absences, European Journal of
Operational Research (2016), http://dx.doi.org/10.1016/j.ejor.2016.04.056}}

\author{Uwe Schauz\\
\small Department of Mathematical Science\\[-0.8ex]
\small Xi'an Jiaotong-Liverpool University\\[-0.8ex]
\small Suzhou 215123, China\\[-0.8ex]
\small \texttt{uwe.schauz@xjtlu.edu.cn}
}

\date{
Mai 20, 2016
}

\usepackage{fixltx2e}

\usepackage[ansinew]{inputenc}
\usepackage{srcltx}
\usepackage[T1]{fontenc}
\usepackage{lmodern}
\usepackage{exscale}

\usepackage{datetime}
\usepackage{ngerman}
\usepackage{textcomp}

\usepackage{tikz}
\usepackage{tkz-graph}
\usepackage{color}
\usepackage{colortbl}
\usepackage{cancel}
\usepackage{float}

\usepackage{amsmath}
\usepackage{amsthm}

\usepackage{amssymb}
\usepackage{stmaryrd}
\usepackage{bm}

\usepackage{enumerate}

\usepackage{afterpage}
\usepackage{xspace}

\newtheoremstyle{Theorem}{.7\baselineskip}{1\baselineskip}{\itshape}{}{\bfseries}{.}{ }{}
\theoremstyle{Theorem}
\newtheorem{Satz}{Theorem}[section]

\newtheorem{Korollar}[Satz]{Corollary}

\newtheorem{Vermutung}[Satz]{Conjecture}

\newtheoremstyle{Definition}{.6\baselineskip}{.8\baselineskip}{}{}{\bfseries}{.}{ }{}
\theoremstyle{Definition}
\newtheorem{Definition}[Satz]{Definition}
\theoremstyle{definition}
\newtheorem{Beispiel}[Satz]{Example}

\theoremstyle{remark}

\newenvironment{sequation*}{\begin{small}\begin{equation*}}{\end{equation*}\end{small}}
\newenvironment{Beweis}[1][Proof]{\begin{proof}[#1]}{\end{proof}}
\newcommand\ps{\small}
\newenvironment{proofsize}[1]{\begin{small}#1}{\end{small}}
\newcommand\psize[1]{\begin{proofsize}#1\end{proofsize}}

\newcommand\textps[1]{\text{\psize{#1}}}

\newenvironment{Liste}[1][]
  {\begin{list}{#1}{\setlength{\topsep}{\itemsep}
   \setlength{\leftmargin}{.5em}\setlength{\labelwidth}{.5em}\setlength{\labelsep}{0em}
   \setlength{\listparindent}{\parindent}\setlength{\parsep}{0ex}}}
  {\end{list}}

  {\setcounter{enumii}{0}\begin{Liste}[\stepcounter{enumii}%
   \emph{#1\arabic{enumi}.\arabic{enumii}}\ssp ,\,\,]}
  {\end{Liste}}

\definecolor{dred}{rgb}{0.9,0,0}

\newcommand\ssp{\kern 1pt}

\newcommand\noms{\hspace{-\mathsurround}}

\newcommand\Rand[1]{
  \marginpar{\raggedleft\scriptsize\hspace{0pt}#1}}%

\renewcommand{\(}{\noms$}
\renewcommand{\)}{\noms$}
\renewcommand\frac[2]{\genfrac{}{}{.4pt}{}{#1}{#2}}
\renewcommand\tfrac[2]{\genfrac{}{}{.4pt}{1}{#1}{#2}}

\newcommand\mathRand[1]{\hspace{\mathsurround}\Rand{#1}\nolinebreak\noms}
\def\rand #1"#2"{\mathRand{\(#2\)}#1#2}
\def\randd #1"#2"#3\randd#4"#5"{\mathRand{\(#2\), \(#5\)}#1#2#3#4#5}

\newcommand\eqby[2][=]%
  {\ensuremath{\overset{\makebox[0pt]{\ensuremath{\smash[t]{\scriptstyle#2}}}}{#1}}}

\renewcommand\o{\emptyset}

\newcommand\Rl{\mathbb{R}}

\newcommand\Z{\mathbb{Z}}
\newcommand\N{\mathbb{N}}

\newcommand\cP{\mathcal{P}}

\renewcommand\sb{\subseteq}
\newcommand\nin{\notin}

\newcommand\sm{\setminus}

\newcommand\mto{\mapsto}
\newcommand\lmto{\longmapsto}
\newcommand\lto{\longrightarrow}

\newcommand\lEqi{\Longleftrightarrow}

\newcommand\ex{\exists\,}

\newcommand\DP{\colon\discretionary{\!\kern -.17em}{}{}}
\newcommand\mitsymbol{\textup{\textbrokenbar}}
\renewcommand\mit{\ssp \ \discretionary{\mitsymbol}{}{}\mitsymbol\ \ssp }

\renewcommand\div{\mathrel{\bigm\lfloor\!\!\!\bigm\lfloor}}
\newcommand\vid{\mathrel{\bigm\rfloor\!\!\!\bigm\rfloor}}
\newcommand\ndiv{\mathrel{\;\!\div\hspace{-12pt}\kern0pt\lower2pt%
  \hbox{\ensuremath{^\diagup}}\!}}
\newcommand\ndivps{\mathrel{\;\!\div\hspace{-9pt}\kern0pt\lower2pt%
  \hbox{\ensuremath{^\diagup}}\!}}

\newcommand\nvid{\mathrel{\;\!\vid\hspace{-12pt}\kern0pt\lower2pt%
  \hbox{\ensuremath{^\diagup}}\!}}
\newcommand\nvidps{\mathrel{\;\!\vid\hspace{-9pt}\kern0pt\lower2pt%
  \hbox{\ensuremath{^\diagup}}\!}}
\providecommand\abs[1]{\lvert#1\rvert}
\providecommand\Abs[1]{\bigl\lvert#1\bigr\rvert}

\DeclareMathOperator\supp{supp}

\newcommand\E{\mathsf{1}}
\newcommand\cop[1]{\lfloor#1\rfloor}

\begin{document}
\maketitle
\begin{abstract}
We study time scheduling problems with allowed absences as a new kind of graph
coloring problem. One may think of a sport tournament where each player (each
team) is permitted a certain number $t$ of absences. We then examine how many
rounds are needed to schedule the whole tournament in the worst case. This upper
limit depends on $t$ and on the structure of the graph $G$ whose edges represent
the games that have to be played, but also on whether or not the absences are
announced before the tournament starts. Therefore, we actually have two upper limits
for the number of required rounds. We have $\chi^t(G)$ for pre-scheduling if all
absences are pre-fixed, and we have $\chi_{\textit{OL}}^t(G)$ for on-line scheduling if
we have to stay flexible and deal with absences when they occur. We conjecture that
$\chi^t(G)=\Delta(G)+2t$ and that $\chi_{\textit{OL}}^t(G)=\chi'(G)+2t.$ The first
conjecture is stronger than the Total Coloring Conjecture while the second is weaker
than the On-Line List Edge Coloring Conjecture. Our conjectures hold for all bipartite
graphs. For complete graphs, we prove them partially. Lower and upper bounds to
$\chi^t(G)$ and $\chi_{\textit{OL}}^t(G)$ for general multigraphs $G$ are established,
too.
\end{abstract}
\newcommand\chio{\chi_{\textit{OL}}}

\section{Introduction}\label{sec.int}

There are many different types of scheduling problems. Some of them arise in pure
mathematics, but many emerge directly out of real-life needs. For example, good
schedules are needed for the assignment of channels or frequencies in communication
networks. They are also needed for the allocation of venues and time slots to the teams
in sport competitions. Most of these problems can be studied as graph coloring
problems, either edge or vertex coloring problems. The graph coloring rule that
adjacent edges (resp.\ vertices) should receive different colors then reflects the most
basic requirement of conflict avoidance, the avoidance of
overlapping appointments in timetables. 
Usually, however, there are additional constraints reflecting additional requirements
and wishes. For instance, in sport league scheduling, one wants to avoid that a team
plays many consecutive games in its hometown. Simultaneously, one wants to
minimize the travel distances of the teams. Moreover, TV networks may want the most
attractive games to be scheduled at certain dates. There is an economic interest
behind many scheduling requirements. Therefore, scheduling has turned into a
research area of its own. Usually, this diverse area is studied in operations research
and computer science. There is a vast literature, see e.g.\ \cite{dk,kkru,lt,rt}\ssp, to
mention but a few. On the web page \cite{kn}\ssp, many references on various topics in
sports scheduling are classified according to different aspects. For mathematical basics
about the theory of graphs and multigraphs (graphs which may have multiple edges
between any two vertices), different coloring concepts and notational foundations, the
reader may consult \cite{di,fiwi,jeto,ya}\ssp.

In the present paper, we examine edge colorings of multigraphs with a new kind of
constraint related to absences. The underlying research should mainly be of interest for
people working in the theory of graph colorings. We hope, however, that our results
and conjectures will also attract interest in the sport scheduling community. In fact, our
research is motivated by time scheduling problems as they arise in sport tournaments
or in the scheduling of timetables at schools. Typically, timetables are set up under the
assumption that everything goes fine and all
participants are available without absences. 
In real life, however, things often do not go as planned. People get sick or otherwise
indisposed. In this case, the best plans can be thrown over. Therefore, it is important to
see how one can deal with absences. Apparently, this problem was not studied in
literature yet, at least not in any systematic way. The closest mathematical concepts, so
far, were \emph{list edge coloring} \cite{bkw,ga,haja}, \emph{on-line list edge coloring}
\cite{schPC,schKp} and \emph{total coloring} \cite{ya}. Our results heavily rely on the
findings in these fields, as we will see. For the general discussion, however, we need
to have mathematical concepts that model time scheduling with absences even closer.
In order to find suitable definitions, we first need to distinguishing two kinds of
absences, \emph{pre-announced pre-fixed absences} and \emph{unannounced
absences}. The following example, with its two parts, illustrates the two types of
absences and their impact on the number of rounds that is needed to accommodate all
games of a tournament. It also explains the
graph-theoretic model that we use
\ssp:

\begin{Beispiel}\label{ex.K3}
\emph{(Part\,1)} Three chess players $A,B,C$ want to play three chess games, the
game $A-B$, the game $A-C$ and the game $B-C$. Each player can play at most one
game per round. Without absences, this can be done in three rounds. One simply has
to play one game per round.

If each player is allowed to miss one round, one round that he has to pre-announce
and prefix, then we may need to arrange one additional round. If player\,\(A\) does not
come to the first round, player\,\(B\) does not come to the second round and
player\,\(C\) does not come to the third round, then three rounds are still enough.
However, in all cases in which at least two players miss simultaneously one of the first
three rounds, a fourth round has to be arranged. The case where player\,\(A\) and
player\,\(B\) cancel the third round and player\,\(C\) cancels the fourth round (whether
or not the fourth round needs to take place) is illustrated in the following assignment of
opponents. Here, a fourth round can actually not be avoided. For the given absences,
four rounds are the minimum, and the presented schedule is optimal in that sense.
Next to the time-table, a corresponding graph coloring is also presented. This
graph-theoretic model shows players as vertices whose color and number indicate the
round in which a player is absent. Games are displayed as edges whose color and
number indicate the round in which the game shall take place\ssp:\smallskip

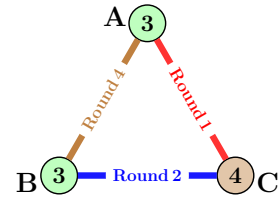
\begin{figure}[H]
\begin{center}
\begin{tabular}{|r||c|c|c|c|}
  \hline
   & \cellcolor{red!25}Round\,1 & \cellcolor{blue!25}Round\,2 & \cellcolor{green!25}Round\,3 & \cellcolor{brown!45}Round\,4 \\\hline\hline
  \cellcolor{gray!7}Player\,\(A\) & \(C\) & free & \textbf{absent} & \(B\) \\\hline
   \cellcolor{gray!7}Player\,\(B\) & free & \(C\) & \textbf{absent} & \(A\) \\\hline
   \cellcolor{gray!7}Player\,\(C\) & \(A\) & \(B\) & free & \textbf{absent} \\\hline
\end{tabular}\hfill
\begin{tikzpicture}
  [transform canvas={scale=1.35}, every node/.style={circle, scale=0.5}]
  \node[draw,fill=green!25,scale=1.1] (a) at (0,1)  {\(\!\bf3\!\)};
  \node[scale=1.3] (a+) at (-0.32,1.07)  {\(\bf A\)};
  \node[draw,fill=green!25,scale=1.1] (b) at (-0.868,-0.5) {\(\!\bf3\!\)};
  \node[scale=1.3] (b+) at (-1.186,-0.57) {\(\bf B\)};
  \node[draw,fill=brown!45,scale=1.1] (c) at (0.868,-0.5) {\(\!\bf4\!\)};
  \node[scale=1.3] (c+) at (1.186,-0.57) {\(\bf C\)};
\SetUpEdge[lw         = 1.8pt,
           labelcolor = white,
           labeltext  = red,
           labelstyle = {font=\sffamily\small,scale=0.8,font=\bf,sloped}]
 \Edge[color=red!80, labeltext=red, label=Round\,1](a)(c)
 \Edge[color=blue!90, labeltext=blue, label=Round\,2](b)(c)
 \Edge[color=brown, labeltext=brown, label=Round\,4](a)(b)
\end{tikzpicture}\hspace{1.8cm}
\caption{Optimal schedule for the indicated pre-fixed absences.}\label{fig.eq1}
\end{center}
\end{figure}\vspace{-8mm}
One can show that four rounds are always enough. For pre-fixed single absences, four
rounds is the upper limit.\smallskip

\emph{(Part\,2)} The situation gets worse if the players do not have to announce their
absences in advance and simply do not show up to one round. In this case, it could
happen that all players come to the first and second round. After these two rounds
there is still at least one game $X-Y$ left over, no mater how the first two rounds are
used. So, two players $X,Y\in\{A,B,C\}$ did not play yet. Now, player\,\(X\) may not
show up to the third round and player\,\(Y\) may not show up to fourth round. For
instance, if the game $A-B$ was not played in the first two rounds ($\{X,Y\}=\{A,B\}$),
this could look as follows:\smallskip

\begin{figure}[H]
\begin{center}
\begin{tabular}{|r||c|c|c|c|}
  \hline
   & \cellcolor{red!25}Round\,1 & \cellcolor{blue!25}Round\,2 & \cellcolor{green!25}Round\,3 & \cellcolor{brown!45}Round\,4 \\\hline\hline
  \cellcolor{gray!10}Player\,\(A\) & \(C\) & free & \textbf{absent} & free \\\hline
  \cellcolor{gray!10}Player\,\(B\) & free & \(C\) & free & \textbf{absent} \\\hline
  \cellcolor{gray!10}Player\,\(C\) & \(A\) & \(B\) & free & free \\\hline
\end{tabular}\hfill
\begin{tikzpicture}
  [transform canvas={scale=1.35}, every node/.style={circle, scale=0.5}]
  \node[draw,fill=green!25,scale=1.1] (a) at (0,1)  {\(\!\bf3\!\)};
  \node[scale=1.3] (a+) at (-0.32,1.07)  {\(\bf A\)};
  \node[draw,fill=brown!45,scale=1.1] (b) at (-0.868,-0.5) {\(\!\bf4\!\)};
  \node[scale=1.3] (b+) at (-1.186,-0.57) {\(\bf B\)};
  \node[draw,fill=gray!10, scale=1.1] (c) at (0.868,-0.5) {\(\!\phantom{\bf4}\!\)};
  \node[scale=1.3] (c+) at (1.186,-0.57) {\(\bf C\)};
\SetUpEdge[lw         = 1.8pt,
           labelcolor = white,
           labeltext  = red,
           labelstyle = {font=\sffamily\small,scale=0.8,font=\bf,sloped}]
 \Edge[color=red!80, labeltext=red, label=Round\,1](a)(c)
 \Edge[color=blue!90, labeltext=blue, label=Round\,2](b)(c)
\end{tikzpicture}\hspace{1.8cm}
\caption{Worst possible unannounced absences after two completed rounds.}\label{fig.eq2}
\end{center}
\end{figure}
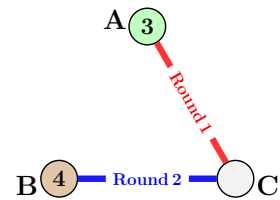\vspace{-8mm}

In this case, a fifth round has to be arranged to accommodate the game $X-Y.$ Since,
at that point, player\,\(X\) and player\,\(Y\) have used up their allowed absences, they
actually will attend the fifth round and the tournament can be concluded there. One can
show that five rounds are always enough. For unannounced single absences, five
rounds is the upper limit.
\end{Beispiel}\smallskip

This example shows that for unannounced absences more rounds might be needed to
accommodate all games, compared to the situation with pre-fixed absences. Out of this
observation, we address two problems. \emph{Pre-scheduling} with pre-announced
pre-fixed absences only, and \emph{on-line scheduling}, where all absences are
unannounced and just happen on the fly. We provide upper and lower bounds on the
number of rounds that is needed to complete all planned games of a tournament. Of
course, this can only be done with some information about the absences. If one team is
absent all the time, we never will finish. Therefore, it seems natural to restrict the
number of absences by some limit. We may permit each player (team) only a certain
number \rand$"t"$ of absences. We may also permit different players $v$ different
numbers $t(v)$ of absences. Apart from that, any choice of matches between the
players is allowed. These matches form the edges of a multigraph \rand$"G".$ With
these notations and parameters, the best general upper bound for the number of
rounds is defined as a new kind of chromatic index. We call this index
\emph{\(t\)-avoiding chromatic index} \randd$"\chi^t"(G),$ respectively \emph{on-line
\(t\)-avoiding chromatic index} \randd$"\chio^t"(G)$ -- one for pre-fixed and one for
unannounced absences. These numbers are the upper limits for the number of rounds
in an optimal scheduling. This means, if we know the number $\chi^t(G),$ resp.\
$\chio^t(G),$ then we know the precise number of required rounds if the absences
appear as unfortunately as possible, within the given frequency limitations $t(v).$

There is also a game-theoretic description of our scheduling problem and the numbers
$\chi^t(G)$ and $\chio^t(G)$. We will not use this approach later on, but we briefly
describe it here, as it clarifies things. The whole scheduling process can be seen as a
meta-game between two meta-players, an \emph{Organizer} and an
\emph{Indisposer}. While Organizer is trying to organize a tournament $G$ within a
certain number $\chi$ of rounds, Indisposer is trying to prevent that by making the
players up to $t$ many times indisposed. There are two versions of that game, one for
pre-announced absences and one for unannounced absences. In the first version,
Indisposer has only one move, in which he determines all absences. He may enter
them into a tabula like the one in Figure\,\ref{fig.eq1}\ssp. Afterwards, Organizer has to
complete the whole schedule in one move by completing the tabula. In the second
version, the only difference is that the tabula is filled column by column. Each round,
Indisposer indicates absences in one column, and then Organizer completes that
column. This could go as in Figure\,\ref{fig.eq2} where $\chi=4$ columns are not
enough to finish the complete tournament $G=K_3$ if $t=1$ many absences are
available in each row. In this game-theoretic setting, the number $\chi^t(G)$, resp.\
$\chio^t(G)$, is the smallest number of columns $\chi$ for which a winning strategy for
Organizer exists.

We can calculate the numbers  $\chi^t(G)$ and $\chio^t(G)$ in several 
cases. In particular, we know $\chi^t(B)$ and $\chio^t(B)$ for all bipartite multigraphs
$B$ and constant or \emph{blockwise constant} $t$ (with $t$ being constant on each of
the two \emph{blocks} of the bipartition). For instance, we may set up a school
timetable and allocate time slots to the lectures that are supposed to be taught.
Assume each teacher may have up to $t_1$ unavailabilities and each class may not be
available at $t_2$ possible lecture times. We can provide the best upper bounds for the
number of time slots that is needed to get all lectures done. In this case, the result only
depends on $t_1$, on $t_2$ and on the maximal degree \rand$"\Delta"(B)$ of the
multigraph $B$ whose edges represent the lectures. Surprisingly, here,
$\chi^t(B)=\chio^t(B).$ This means that upfront notification of absences does not give
us an advantage, not for the general upper bound $\chi^t(B).$ This is an exceptional
phenomenon in the bipartite case.

For general graphs $G$, we only have partial results. Sometimes, the upper bound
$\chio^t(G)$ can be improved if absences are reported before the whole tournament
starts. For instance, this seems to be the case for all complete graphs \rand$"K_n"$
with an odd number of vertices. Actually, complete graphs are of special interest for
round-robin tournaments. Therefore, some attention is given to complete graphs. We
can provide partial result here. If each player is allowed to miss just one announced
round, we can provide detailed information, but we do not know $\chi^2(K_n)$ for odd
$n,$ for example. If the absences are unannounced, we can calculate $\chio^t(K_n)$
for all constant multiplicities $t,$ but only if $n-1$ is even or prime. Our results suggest
two conjectures about $\chi^t(G)$ and $\chio^t(G)$ for general graphs $G$ if $t$ is
constant, Conjecture\,\ref{con.ub} and Conjecture\,\ref{con.ubo}\ssp. The first one,
about $\chi^t(G),$ is a reasonable strengthening of the Total Coloring Conjecture. It
says that $\chi^t(G)=\Delta(G)+2t.$ If correct, it would give us $\chi^t(G)$ for all
\emph{simple} graphs (graphs without multiple edges). The second one, about
$\chio^t(G),$ is a weakening of the on-line version of the List Edge Coloring
Conjecture. It says that $\chio^t(G)=\chi'(G)+2t.$ This means, if the on-line
strengthening of the long-standing List Edge Coloring Conjecture should hold, then we
know $\chio^t(G)$ for all multigraphs $G$ (\ssp at least if we know the chromatic index
\rand$"\chi'"(G)$).

Following this introduction, we present two structurally similar sections. In
Section\,\ref{sec.CI}\ssp, we study the situation with pre-announced pre-fixed
absences. In Section\,\ref{sec.OL}\ssp, we study the on-line situation with
unannounced absences. In both sections, we start with the basic definitions, then look
at lower bounds (Section\,\ref{sec.CIlb}\ssp, resp.\ \ref{sec.OLlb}), upper bounds
(Section\,\ref{sec.CIub}\ssp, resp.\ \ref{sec.OLub}), results for bipartite graphs
(Section\,\ref{sec.CIB}\ssp, resp.\ \ref{sec.OLB}) and, finally, results for complete
graphs (Section\,\ref{sec.CIKn}\ssp, resp.\ \ref{sec.OLKn}). The results in the two
sections about bipartite graphs $B$ are very much alike. We explain the basics about
$\chi^t(B)$ in Section\,\ref{sec.CIB} and keep this section free from the on-line index
$\chio^t(B).$ The additional background for the study of $\chio^t(B)$ is then introduced
in Section\,\ref{sec.OLB}\ssp, without repetition of the basics.
To keep our presentation systematic, however, we state the main results about bipartite
graphs there again, just generalized to $\chio^t(B).$ The situation, however, is very
different for complete graphs. The results in the two sections about complete graphs
are rather unlike. Therefore, we compare the main differences
between pre-scheduling and on-line scheduling in 
Section\,\ref{sec.comp}\ssp. For complete graphs and $t=1$, the differences between
$\chi^1(K_n)$ and $\chio^1(K_n)$ are discussed in Example\,\ref{ex.Kn} and
Figure\,\ref{fig.2}\ssp, which extend the basic Example\,\ref{ex.K3} in this introduction.
Finally, in Section\,\ref{sec.conc}\ssp, we summarize the provided evidence for our two
conjectures and point out some interesting aspects and open cases.

\section{Pre-Fixed Absences and Pre-Scheduling}\label{sec.CI}

We model time scheduling problems between competitors in sport tournaments using
graph theory\ssp: In a multigraph \rand$"G=(V,E)"$ the vertex set $V$ represents the
set of competitors or competing teams. An edge \rand$"uv":=\{u,v\}$ of the edge
multiset $E$ stand for a match between $u$ and $v.$ Single absences of players may
be described by a vertex labeling \rand$"c"\DP V\!\lto\Z^+$ with positive integers --
the label $c(v)$ of a player $v$ representing the round in which that player is not
available. For example, $c(v)=5$ would mean that $v$ is not available in the
$5^\text{th}$ round, if the tournament should last $5$ rounds. The round in which a
game $e\in E$ shall be played is denoted by \rand$"c'"(e).$ Hence, we are looking for
edge labelings $c'\DP E\lto\Z^+$ that go well together with the absences given by $c.$
We say that an edge labeling $c'\DP E\lto\Z^+$ \emph{avoids} $c$ if the \emph{color}
$c'(uv)$ of every edge $uv\in E$ differs from the labels $c(u)$ and
$c(v)$ of its ends. In other words, for all $e\in E,$%
 \rand\begin{equation}
c'(e)\,\nin\,"c(e)":=\{\ssp c(v)\mit v\in e\ssp \}.
\end{equation}
More generally, one may allow multi-labelings \rand$c\DP V\!\lto"\cP(\Z^+)":=\{I\!\mit\!
I\sb\Z^+\}$ to describe multiple absences. Then, the set of players that do not show up
in round number $i$ is \rand\begin{equation}\label{eq.Ui} "U_i"\,:=\,\{\ssp v\in V\mit i\in
c(v)\ssp \},
\end{equation}
provided that the tournament is not already finished before the $i^\text{th}$ round. In
this situation, \emph{\(c\)-avoidance} means that, for all $e\in E,$
 \rand\begin{equation}
c'(e)\,\nin\,"c(e)":=\bigcup_{v\in e}c(v).
\end{equation}
Avoiding $c$, however, is not the only requirement and objective. Since a player $v$
cannot play against two other players at the same time, we need to find (proper)
\emph{edge colorings} $c'\DP E\lto\Z^+\!,$ i.e., we need that $c'(e)\neq c'(f)$ for every
two edges $e$ and $f$ with $e\cap f\neq\o.$ We define the \emph{\(c\)-avoiding
chromatic index} \rand$"\chi'_c"(G)$ of $G$
 as the smallest number of rounds that is needed to arrange all matches
$e\in E$ in a way that avoids the absences given by $c.$ Mathematically,
\begin{equation}\label{eq.xc}
\chi'_c(G)\,:=\,\min\{\ssp m\in\Z^+\mit\textps{$\ex$ edge coloring $c'\DP E\to(m]$ {\ps that avoids} $c$}\ssp \},
\end{equation}
where
 \rand\begin{equation}
"(m]"\,=\,(0,m]\,:=\,\{1,2,\dotsc,m\}.
\end{equation}
If $\chi'_c(G)\leq m,$ we say that $G$ is \emph{\(c\)-avoiding edge \(m\)-colorable}. In
other words, $G$ is \(c\)-avoiding edge \(m\)-colorable if and only if $G$ is \emph{list
edge colorable} with respect to the color lists $L_e$ given by
 \rand\begin{equation}\label{eq.Le} "L_e"\,:=\,(m]\!\sm\!c(e).
\end{equation}
It means that it is possible to find a proper edge coloring by choosing the colors of the
edges $e\in E$ from the corresponding lists $L_e.$ In short, there exists an edge
coloring $c'$ of the form
\begin{equation}
c'\DP E\ni e\lmto c'(e)\in L_e.
\end{equation}
This observation will become helpful later. It will allow us to utilize list coloring
theorems. However, it would be very surprising if $\chi'_c(G)$ could be expressed by
any general formula of reasonable complexity. It will usually depend on $c$ in
complicated ways, as we will already see in our study of complete graphs. Therefore,
the main purpose of this paper is to provide upper bounds for $\chi'_c(G),$ if some
partial information about $c$ is given. Here, a quite reasonable restriction on the
absences given by $c$ is that every player $v$ can miss only a limited number
\rand$"t"$ of rounds. To investigate this situation, we call a multi-labeling $c\DP
V\!\to\cP(\Z^+)$ a \emph{\(t\)-labeling} if $\abs{c(v)}\leq t$ for all $v\in V\!.$ We are
interested in the \emph{\(t\)-avoiding chromatic index}
 \rand\begin{equation}\label{eq.xt}
"\chi^t"(G)\,:=\,\max\{\ssp \chi'_c(G)\mit\textps{$c$ is a \(t\)-labeling}\ssp \}.
\end{equation}
For $t=0,$ this is just the \emph{chromatic index}
 \rand\begin{equation}
\chi^0(G)\,=\,"\chi'"(G)\,:=\,\min\{\ssp m\in\Z^+\mit\textps{$\ex$ edge coloring $c'\DP
E\to(m]$}\ssp \}.
\end{equation}
For $t=1,$ it can be written as
\begin{equation}\label{eq.x1}
\chi^1(G)\,=\,\max\{\ssp \chi'_c(G)\mit c\DP V\!\to\Z^+\ssp \}.
\end{equation}
If $\chi^t(G)\,\leq\,m,$ we say that $G$ is \emph{\(t\)-avoiding edge \(m\)-colorable}. It is
also possible to generalize this concept to \emph{multiplicity functions}
\rand$"t"\DP V\!\to\N:=\{0,1,\dotsc\}.$ 
In this case, however, note that the term ``\(t\)-avoiding'' has a different meaning from
the previously defined term ``\(c\)-avoiding''; even though $t=c$ is possible, as
multiplicity functions $t$ are also vertex labelings. This is the reason why we write $t$
in the exponent of $\chi.$ We need to distinguish $\chi^t$ from $\chi'_c.$

\subsection{General Lower Bounds (pre-scheduling)}\label{sec.CIlb}

Since, in the definition of $\chi^t(G)$ in Equation\,\eqref{eq.xt}\ssp, one also has to
consider the multi-labeling $c$ with $c(v):=(t]=\{1,2,\dotsc,t\}$ for all $v\in V,$ we see
that, for non-edgeless graphs $G,$
\begin{equation}\label{eq.lb2x}
\chi^t(G)\,\geq\,\chi'(G)+t.
\end{equation}
Usually, however, this is far from being the best lower bound. If we assign the set $(t]$
to a
vertex $v_0$ of maximal degree, \rand\rand$"d"(v_0)="\Delta"(G)
,$ and assign the set $(2t]\!\sm\!(t]$ to all its neighbors $u,$ then we need
$\Delta(G)+2t$ many colors for the connecting edges $v_0u$ alone. Hence, for
non-edgeless graphs $G,$
\begin{equation}\label{eq.lb2}
\chi^t(G)\,\geq\,\Delta(G)+2t.
\end{equation}
For simple graphs and $t\geq1,$ this lower bound is better than our initial lower bound
$\chi'(G)+t,$ because $\Delta(G)\geq\chi'(G)-1$ for simple graphs, by Vizing's
Theorem. Moreover, it can be generalized to non-constant $t.$ Obviously, in that case,
\begin{equation}\label{eq.lb3}
\chi^t(G)\,\geq\,\Delta(G)+\min_{uv\in E}(t(u)+t(v)),
\end{equation}
but note that all our lower bounds do not hold for edge-less graphs, as then
$\chi^t(G)=0.$

If, in the definition of $\chi^1(G)$ in Equation\,\eqref{eq.x1}\ssp, one restricts $c\DP
V\!\to\Z^+$ to (proper) vertex colorings, one almost obtains the so-called \emph{total
chromatic number} \cite{ya}\ssp, which is defined by
 \rand\begin{equation}\label{eq.xT}
"\chi_T"(G)\,:=\,\min\{\ssp m\in\Z^+\mit\text{$\ex c\in(m]^V\!,\,c'\in(m]^E$ {\ps both
proper,} $c'$ {\ps avoids} $c$}\ssp \}.
\end{equation}
The only reason why $\chi_T(G)$ could be bigger than $\chi^1(G)$ is that the $m$ in
its definition must be big enough to allow a vertex coloring $c$ of $G,$ which lifts
$\chi_T(G)$ up to $\chi(G)$ at the outset. The requirement $m\geq\chi(G)$ is not
explicitly requested in Equation\,\eqref{eq.xc}\ssp, the definition of $\chi'_c(G)$
underlying $\chi^1(G).$ 
Using Lower Bound\,\eqref{eq.lb2x} and Brooks' Theorem\ssp, however, we actually
can deduce that $\chi^1(G)\geq\chi'(G)+1\geq\Delta(G)+1\geq\chi(G),$ for
non-edgeless graphs $G.$ So, the initial lift of $\chi_T(G)$ up to $\chi(G)$ does not lift
it above $\chi^1(G).$ Hence, for non-edgeless graphs $G,$ we always have
\begin{equation}\label{eq.lbT1}
\chi^1(G)\,\geq\,\chi_{_T}(G).
\end{equation}
The two parameters are not equal though, in general, e.g.\
$\chi^1(K_3)=4>3=\chi_{_T}(K_3).$

\subsection{General Upper Bounds (pre-scheduling)}\label{sec.CIub}

As already recognizable in Equation\,\eqref{eq.Le} above, an interesting connection is
given to list colorings. Therefore, we consider the \emph{list chromatic index}. It is
defined as
 \rand\begin{equation}\label{eq.xl}
"\chi_\ell'"(G)\,:=\,\min\{\ssp m\in\Z^+\mit G\textps{\ is \(m\)-list edge colorable}\ssp \},
\end{equation}
where $G$ is \emph{\(m\)-list edge colorable} if it is edge colorable from every
\emph{\(m\)-lists system} \rand$e\mto "L_e",$ i.e.\ every system of \emph{lists} (sets)
$L_e$ with $\abs{L_e}=m$ available \emph{colors} at each edge $e\in E.$ Since, for
\(t\)-labelings $c$ and all $e\in E,$ $\abs{c(e)}\leq 2t$ and
\begin{equation}
\Abs{(m]\!\sm\!c(e)}\,\geq\,m-2t,
\end{equation}
we see that $G$ is \(t\)-avoiding edge \(m\)-colorable if $G$ is \((m{-}2t)\)-list edge
colorable. Hence, always
\begin{equation}\label{eq.UBLC1}
\chi^t(G)\,\leq\,\chi_\ell'(G)+2t.
\end{equation}
If the List Edge Coloring Conjecture holds for $G,$ i.e.\ if $\chi_\ell'(G)=\chi'(G),$ then
\begin{equation}\label{eq.UBLC}
\chi^t(G)\,\leq\,\chi'(G)+2t.
\end{equation}
For bipartite multigraphs, this actually holds, as we will discuss in the next section. In
this case, we may also replace the chromatic index $\chi'$ with the maximal degree
$\Delta$ in the upper bound, as then $\chi'(G)=\Delta(G)$ by König's Theorem.
Therefore, one may wonder if this is possible for other types of graphs as well. As we
will see, if we have multiple edges in a non-bipartite graph $G,$ an upper bound of
$\Delta(G)+2t$ does not hold in general. Of course, it does not hold for simple graphs if
$t=0$ either, because then Vizing's upper bound
$\Delta(G)+1\geq
\chi'(G)=\chi^0(G)$ is best possible; e.g.\ $\chi'(K_3)=\Delta(K_3)+1.$ 
For $t=1$ and simple graphs $G,$ however, the upper bound $\Delta(G)+2$ appears
as a reasonable strengthening of \emph{\null\!Vizing's and Behzad's Total Coloring
Conjecture}, which says that
\begin{equation}
\chi_{_T}(G)\,\leq\,\Delta(G)+2.
\end{equation}
We will prove this strengthening for complete graphs. As explained, it also holds for
bipartite graphs. If it should be accurate for all simple graphs $G$ then possibly even
$\chi^t(G)\leq\Delta(G)+2t$ for arbitrary $t>0.$ Together with Lower
Bound\,\eqref{eq.lb2} this would mean that $\chi^t(G)=\Delta(G)+2t,$ so that
$\chi^t(G)$ could be easily determined for all graphs without multiple edges. Since this
would be a surprisingly concise result, we want to encourage further research and
venture the following conjecture\ssp :

\begin{Vermutung}\label{con.ub}
For every non-edgeless simple graphs $G$ and every $t\in\Z^+\!\,$
$$
\chi^t(G)\,=\,\Delta(G)+2t.
$$
\end{Vermutung}

However, the best general result that we can provide here arises from Shannon’s
Bound for arbitrary multigraphs, or from Borodin, Kostochka and Woodall's sharpening
\cite[Theorem\,4\ssp ]{bkw} of that bound. Their sharpening says that multigraphs are
list edge colorable if, for all $uv\in E,$ the list of $uv$ contains at least
$\max\bigl(d(u),d(v)\bigr)+\cop{\tfrac{1}{2}\min\bigl(d(u),d(w)\bigr)}$ many colors, where
\rand$"\cop{x}"$ denotes the biggest integer below $x\in\Rl.$ This leads to the
following result\ssp :

\begin{Satz}\label{sz.sh}
For every multigraph $G=(V,E)$ and every function $t\DP V\!\to\N,$
$$
\chi^t(G)\,\leq\,\max_{uv\in E}\bigl(\max(d(u),d(v))+\cop{\tfrac{1}{2}\min(d(u),d(v))}+t(u)+t(v)\bigr).
$$
In particular, if $t$ is constant,
$$
\chi^t(G)\,\leq\,\cop{\tfrac{3}{2}\Delta(G)}+2t.
$$
\end{Satz}

Here, for multigraphs, it is not hard to see that the factor $\tfrac{3}{2}$ in the upper
bound $\cop{\tfrac{3}{2}\Delta(G)}+2t$ cannot be improved. We have the following
obvious example\ssp:
\begin{Beispiel}
Consider the \emph{thick triangle} $G:=s\times K_3$ (with $s$ parallel edges between
any two of its 3 vertices $v_1,v_2,v_3$), and take the three sets
$c(v_1):=\{1,2,\dotsc,2r\},$ $c(v_2):=\{1,2,\dotsc r\,,\,2r{+}1,2r{+}2,\dotsc,3r\}$ and
$c(v_3):=\{r{+}1,r{+}2,\dotsc,3r\}$ as sets of absences (\ssp i.e.\ $t=2r$), then
$$
\chi'_c(G)\,=\,3s+3r\,=\,\tfrac{3}{2}\Delta(G)+\tfrac{3}{2}t.
$$
\end{Beispiel}

\subsection{Bipartite Tournaments (pre-scheduling)}\label{sec.CIB}

Galvin could show in \cite{ga} that the List Edge Coloring Conjecture holds for bipartite
multigraphs \rand$"B"\,$ and
\begin{equation}
\chi_\ell'(B)\,=\,\chi'(B)\,=\,\Delta(B).
\end{equation}
Therefore, Upper Bound \eqref{eq.UBLC1} coincides with Lower
Bound\,\eqref{eq.lb2}\ssp, and we obtain the following Theorem\ssp :

\begin{Satz}\label{sz.ga}
For every non-edgeless bipartite multigraph $B$ and every $t\in\N,$
$$
\chi^t(B)\,=\,\chi'(B)+2t\,=\,\Delta(B)+2t.
$$
\end{Satz}

Moreover, one may even use Borodin, Kostochka and Woodall's sharpened version
\cite[Theorem\,3\ssp ]{bkw} of Galvin's Theorem, which says that $B$ is list edge
colorable if the list of every edge $uv$ contains at least $\max\bigl(d(u),d(v)\bigr)$
many colors. This yields the following improvement of Theorem\,\ref{sz.ga}\ssp :

\begin{Satz}\label{sz.BKW}
For every bipartite multigraph $B=(V,E)$ and every function $t\DP V\!\to\N,$
$$
\chi^t(B)\,\leq\,\max_{uv\in E}\bigl(\max(d(u),d(v))+t(u)+t(v)\bigr)\,\leq\,\Delta(B)+\max_{uv\in E}\bigl(t(u)+t(v)\bigr).
$$
\end{Satz}

In this theorem, the expression in the middle depends on the number of allowed
absences $t(u)$ and the ``busyness'' $d(u)$ of the players $u$ and their neighbors.
Hence, to keep the number of required rounds $\chi^t(B)$ small, it could be good to
restrict the number of absences of the busiest players more strongly. 
In real-life situations, however, it might be more realistic to expect that blockwise
constant functions $t$ are used. Indeed, in the previously mentioned example with
teachers and school classes, this is a very reasonable scenario. All teachers have the
same number $t_1$ of allowed absences, and all classes have the same number
$t_2$ of allowed absences. If in this situation the teachers are not available during the
first $t_1$ rounds, and the classes are absent during the following $t_2$ rounds, then
$t_1+t_2$ rounds are lost. Hence, if the underlying graph $B$ is non-edgeless, then
$\chi'(B)+t_1+t_2$ rounds are required, i.e., in accordance with Lower
Bound\,\eqref{eq.lb3}\ssp,
\begin{equation}\label{eq.lbC2}
\chi^t(B)\,\geq\,\chi'(B)+t_1+t_2\,=\,\Delta(B)+t_1+t_2.
\end{equation}
This coincides with the given upper bound. We obtain the following corollary\ssp :

\begin{Korollar}\label{cor.BKW}
If $B$ is non-edgeless and $t\DP V\!\to\N$ takes the value $t_1$ on one block and the
value $t_2$ on the other block of the bipartition, then
\begin{equation}\label{eq.lbC3}
\chi^t(B)\,=\,\chi'(B)+t_1+t_2\,=\,\Delta(B)+t_1+t_2.
\end{equation}
\end{Korollar}

\subsection{Complete Tournaments (pre-scheduling)}\label{sec.CIKn}

The results about list edge colorings that we used in the study of bipartite multigraphs
are quite strong, compared to what we know for other graphs. Actually, at present, it is
not even known if the List Edge Coloring Conjecture holds for all complete graphs.
H\"{a}ggkvist and Janssen \cite{haja} provided the upper bound
\begin{equation}
\chi'_\ell(K_n)\,\leq\,n,
\end{equation}
for all $n\in\Z^+\!\!$, which implies that
\begin{equation}
n+2t-1\,\leq\,\chi^t(K_n)\,\leq\,n+2t.
\end{equation}
Here, the upper bound is already close to the given lower bound, which arises from
Lower Bound\,\eqref{eq.lb2}\ssp. 
For all even $n$ for which $\chi'_\ell(K_n)=\chi'(K_n),$ however, the upper bound can
be improved by one, because $\chi'(K_n)=n-1$ for even $n.$ The problem is that the
conjectured equality $\chi'_\ell(K_n)=\chi'(K_n)$ is hard to prove in general. In
\cite{schKp}\ssp, we could only prove it under the additional assumption that $n-1$ is
an odd prime. Moreover, from the outset, this approach to improve the upper bound
$n+2t$ is restricted to even $n,$ because $\chi'(K_n)=n$ for odd $n.$ Nevertheless,
the upper bound might be improvable for all $n,$ as Conjecture\,\ref{con.ub} suggests.
Indeed, at least for $t=1,$ this is the case\ssp :

\begin{Satz}\label{sz.x1Kn}
For all $n\geq2,$
$$\chi^1(K_n)\,=\,n+1.$$
\end{Satz}

\begin{figure}[t]
\begin{center}
\includegraphics[scale = .27
]{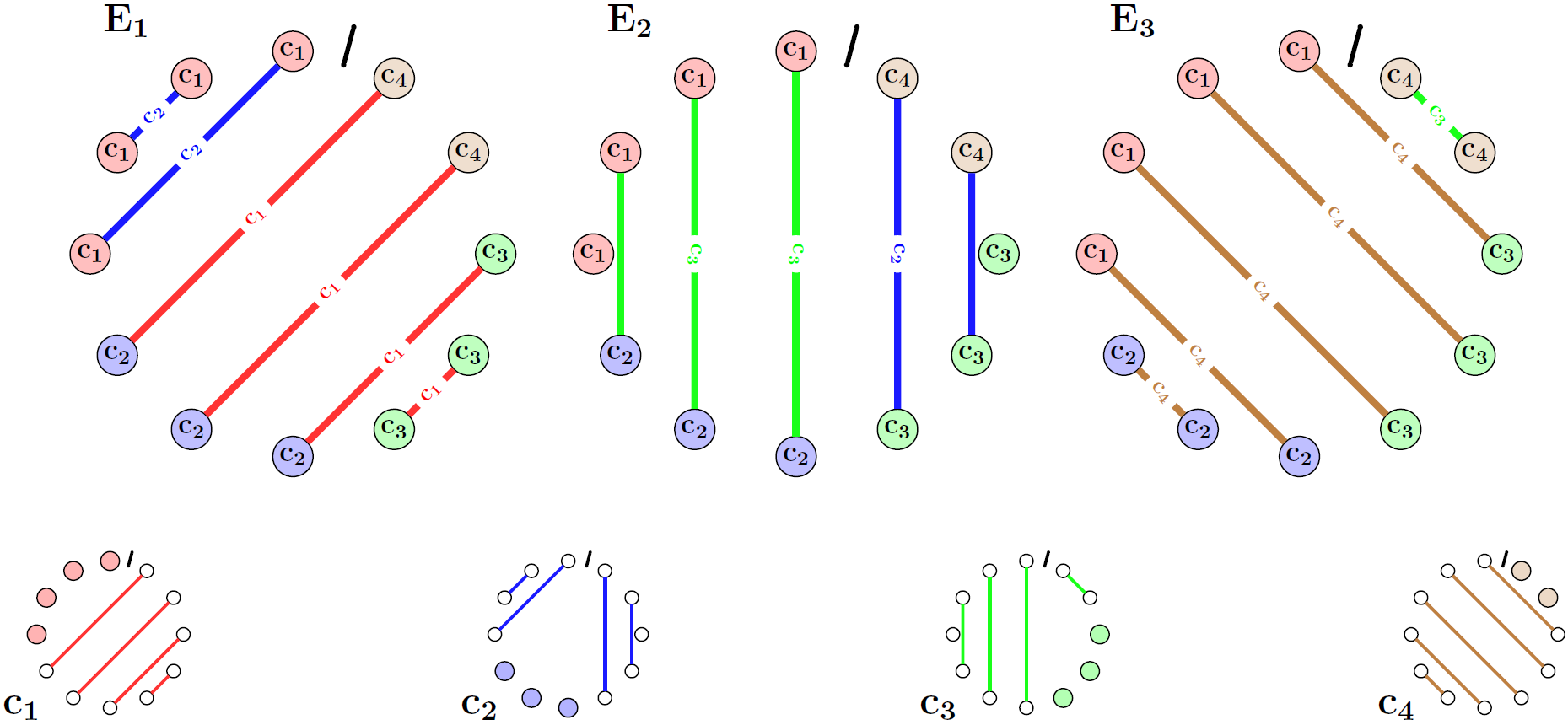}
\caption{The $s-1$ bichromatic parallel classes $E_j$ and the $s$ colors $c_j$ of
the vertices.}\label{fig.1}
\end{center}\vspace{-3ex}
\end{figure}

\begin{Beweis}
Based on Lower Bound\,\eqref{eq.lb2} it suffices to show that $\chi^1(K_n)\leq n+1.$
This holds for $n=2.$ Thus, assume $n\geq3$ and let $c\DP(n]\to\Z^+$ be a labeling of
$V(K_n)=(n].$ We show that $G$ is \(c\)-avoiding edge \((n+1)\)-colorable. Without loss
of generality, we assume that $c((n])\sb(n{+}1].$ Let $c_1,c_2,\dotsc,c_s$ denote the
different values of $c,$ i.e.\ $c((n])=\{c_1,c_2,\dotsc,c_s\}.$ After permuting the vertices,
we may suppose that $c_1$ occurs on the first $n_1$ vertices, $c_2$ on the next
$n_2,$ etc. For $j=1,2,\dotsc,s,$ set $m_j:=n_1+n_2+\dotsb+n_j.$

To determine the color classes of a suitable coloring, we arrange the $n$ vertices
equidistantly and counter-clockwise around a cycle, and draw the edges as straight
lines. Then, the edges form $n$ classes of parallel edges. We denote the parallel class
of an edge $e$ as $[e].$ For instance, in $K_6,$ $[\{1,2\}]=\{\{1,2\},\{3,6\},\{4,5\}\}$ and
$[\{1,3\}]=\{\{1,3\},\{4,6\}\}.$ We say that a subset $f\sb(n]$ is bigger than a subset
$e\sb(n]$ if $\min(f)>\max(e),$ e.g. the edge $\{3,6\}$ is bigger than the edge $\{2,1\}.$

The initial idea is to color the parallel classes monochromatically, using one exclusive
color for all the edges inside one class only. In this way, we obtain a proper edge
coloring of $K_n$ with $n$ colors. However, this coloring usually does not avoid $c.$
We only use it in the case $s=1.$ If $s>1,$ we modify our approach. We color each of
the $s{-}1$ parallel classes $E_1:=[\{1,m_1\}]$, \dots, $E_{s{-}1}:=[\{1,m_{s{-}1}\}]$ with
two colors (where we set $E_1:=[\{2,n\}]$ if $m_1=1$). We assign the color $c_j$ to all
edges in $E_j$ that are bigger than the set $\{1,m_j\},$ and the color $c_{j+1}$ to the
remaining edges in $E_j$.  This is shown in an example in Figure\,\ref{fig.1}, where we
have $s-1=3$ many bichromatic parallel classes. Using this coloring strategy, each
color $c_j$ with $1<j<s$ is used
only in $E_j$ and $E_{j-1}.$ 
The \(c_j\)-colored edges in $E_{j-1}$, however, connect only vertices in $(m_{j-1}],$
while the \(c_j\)-colored edges in $E_j$ connect only vertices in the complement of
$(m_j].$ Hence, there is no color conflict between these edges. We also avoid the
vertex labeling $c,$ which shows the label $c_j$ only inside $(m_j]\!\sm\!(m_{j-1}].$
There are no problematic incidences or adjacencies, also not within the colors $c_1$
and $c_s\ssp.$ This is illustrated in the example in Figure\,\ref{fig.1} (were we have
$s=4$ many color classes, displayed in the four smaller sub-diagrams). Now, the
colors in $(n{+}1]\sm\{c_1,c_2,\dotsc,c_s\}$ can be used to color the remaining
$n-(s-1)$ many parallel classes monochromatically. The resulting proper edge coloring
avoids $c.$ So, $\chi'_c(K_n)\leq n+1$ and $\chi^1(K_n)=n+1.$
\end{Beweis}

This theorem can also be reinterpreted in terms of Latin squares. In the study of the
case $t=1,$ we only have to consider \(1\)-labelings $c\DP V(K_n)\to\Z^+\!.$ Now, if
$V(K_n)=(n],$ the labels are $c(1),c(2),\dotsc,c(n),$ and we may write them in the
diagonal of an $n\times n$ matrix. The colors $c'(ij)$ of the edges $ij$ in a
\(c\)-avoiding coloring $c'$ of $K_n$ can then be written in the positions $(i,j)$ and
$(j,i)$ of that matrix. In this way we get a symmetric partial Latin square with $c$ in the
diagonal, i.e.\ no column or row contains an entry twice. It is called \emph{partial},
because it may use more than $n$ different entries, and may then be viewed as part of
a possibly existent bigger Latin square. In this matrix interpretation,
Theorem\,\ref{sz.x1Kn} says the following\ssp :

\begin{Korollar}
For every function $c\DP(n]\to(n{+}1]$, there exists a symmetric partial Latin square in
$(n{+}1]^{n\times n}$ with the sequence $c(1),c(2),\dotsc,c(n)$ in the diagonal.
\end{Korollar}

This Corollary seems to be new, but 
there is a known improvement if $c$ has certain properties\ssp :

\begin{Satz}\label{sz.sl}
Let $c\DP(n]\to(n]$ be a function that takes $s$ many different values
$c_1,c_2,\dotsc,c_s,$ and let $n_1,n_2,\dotsc,n_s$ be the frequencies of these values.
Then there exists a symmetric Latin square in $(n]^{n\times n}$ with the sequence
$c(1),c(2),\dotsc,c(n)$ in the diagonal if and only if $n_1\equiv n_2\equiv\dotsb\equiv
n_s\equiv n\pmod{2}.$
\end{Satz}

\begin{Beweis}
This is a special case of \cite[Theorem\,1]{cr} and of \cite[Corollary\,12]{an}\ssp.
\end{Beweis}

We can use this theorem to say more about the number of rounds that is needed to
schedule a complete tournament with at most one absence per player\ssp :

\begin{Satz}\label{sz.xcKn}
Let  $c\DP(n]\to\Z^+$ be a labeling of the vertex set $(n]$ of $K_n,$ $n\geq2.$ Assume
that $c$ takes $s$ many different values $c_1,c_2,\dotsb,c_s$ inside $(n]\sb\Z^+\!.$
Let $n_1,n_2,\dotsc,n_s$ be the frequencies of these values, and let
$n^+:=n-(n_1+n_2+\dotsb+n_s)$. Then
$$
\chi'_c(K_n)=\!\begin{cases}
 n-1&\!\!\text{if \,$n $\, is even and $c(j)\geq n$ for all $j\in(n],$}\\
 n&\!\!\text{else if \,$n_i\equiv n\!\!\!\pmod{2}$\, for at least \,$s-n^+\!$ many \,$i\in(s],$}\\
 n+1&\!\!\text{else.}
\end{cases}
$$
\end{Satz}

\begin{Beweis}
By Theorem\,\ref{sz.x1Kn}\ssp, we already know that
\begin{equation}
\chi'_c(K_n)\,\leq\,n+1.
\end{equation}
Theorem\,\ref{sz.sl} can be used to decide precisely when $\chi'_c(K_n)\leq n:$

The interpretation of Theorem\,\ref{sz.sl} as a statement about \(c\)-avoiding edge
colorings of $K_n$ is just opposite to what we have done above. If $c$ is a map into
$(n],$ the condition of Theorem\,\ref{sz.sl} is necessary and sufficient for
$\chi'_c(K_n)\leq n.$ If $c$ takes values bigger than $n$ (i.e.\ $n^+>0$), then one may
change these overstepping function values into suitable values inside $(n]$, in order to
be able to apply Theorem\,\ref{sz.sl}\ssp. If this strategy succeeds and only $n$ rounds
are required eventually, then it does not matter that we changed these values. We
reach the end of the game before the original overstepping values come into effect.
Therefore, we only have to see how we can choose suitable replacement values in this
reassignment strategy. To start with, assume that the frequency \randd$"n_{i_0}"$ of
just one value \randd$"c_{i_0}"\leq n$ has the wrong parity, $n_{i_0}\not\equiv
n\pmod{2}.$ Now, if $n$ is even, then the number $n^+\!$ of all $j\in(n]$ with $c(j)>n$
is necessarily odd, under our assumption. In this case, one may change the $n^+\!$
overstepping function values into $c_{i_0},$ to increase the frequency $n_{i_0}$ to
$n_{i_0}+n^+\!.$ If $n$ is odd then $n^+\!$ can be both, odd or even. If $n^+\!$ is odd,
then we proceed as before. If $n^+\!$ is even, then we split $n^+\!$ into two odd
summands $n^+_1$ and $n^+_2.$ We increase $n_{i_0}$ to $n_{i_0}+n^+_1$ as
before, and introduce a new function value in $(n]$ with frequency $n^+_2.$ In all
cases, the new frequencies have the right parity, and Theorem\,\ref{sz.sl} guarantees
that $\chi'_c(K_n)\leq n.$ It is not hard to see that for $k$ many wrong parities, $k$
many $j\in(n]$ with $c(j)>n$ suffice. Put another way, having at least $s-n^+\!$ many
times the correct parity suffices. From the necessity of the parity condition in
Theorem\,\ref{sz.sl}\ssp, it is also not hard to see that this sufficient condition is also
necessary. Indeed, if $n$ colors are used for the edges, then every vertex will be
incident to exactly $n-1$ colors. So, the one missing color may be placed in the
diagonal, if it was not already there, as a pre-fixed absence inside $(n]$. Hence,
overstepping function values necessarily have a ``fortunate'' reassignment possibility
inside $(n]$, if a \(c\)-avoiding edge \(n\)-coloring exists at all. This shows that
$\chi'_c(K_n)\leq n$ if and only if $n_i\equiv n\pmod{2}$ for at least $s-n^+$ many
$i\in(s].$ It remains to study the cases with $\chi'_c(K_n)\leq n-1:$

Actually, if $n$ is odd, this case cannot occur, as then $\chi'_c(K_n)\geq\chi'(K_n)=n$. 
If $n$ is even, then $\chi'(K_n)=n-1$, and $\chi_c'(K_n)\geq n-1$ follows. We have
$\chi_c'(K_n)=n-1$, if $n$ is even and if $c$ avoids conflicts with the edge colors of
any one edge \((n-1)\)-coloring. However, as in every edge \((n{-}1)\)-coloring all $n-1$
colors occur at the edges of every vertex of $K_n$, this means that $c(j)\geq n$ for all
$j\in(n],$ as claimed in the theorem.
\end{Beweis}


\section{Unannounced Absences and On-Line Scheduling}\label{sec.OL}

In this section, we investigate what happens when we do not know in advance in which
rounds the players actually appear. In real life, the organizers may insist on a certain
warning time, to be able to prepare the next round in time. In our investigations,
however, we assume that we find out who attends right before the start of each round.
There is still a given upper bound $t$ for the number of allowed absences, possibly a
function $t\DP V\!\to\N,$ but the precise allocation of the absences is not fixed
beforehand.

We say that a tournament is \emph{on-line \(t\)-avoiding edge \(m\)-colorable} if it is
always possible to schedule the matches in such a way that we finish the whole
tournament after at most $m$ rounds. This is a simple definition in the language of
sport schedulings, but it is a bit harder to translate this notion into mathematical terms.
We observe that a tournament $G=(V,E)$ can be scheduled in $m$ rounds if and only
if the first round can be scheduled in such a way that the remaining games can be
scheduled in $m-1$ rounds, no matter which subset \rand$"U"$ of players is absent in
the first round. More precisely, there must exist a set of \emph{independent} (pairwise
disjoint) matches \rand$"F"$ in $G-U$ such that the games in $G-F$ can be scheduled
in $m-1$ rounds, for every set $U$ of players who are allowed to miss a round, i.e.\ for
all
 \rand\begin{equation}
U\,\sb\,"\supp"(t)\,:=\,\{\ssp v\in
      V\mit t(v)\neq0\ssp \}.
\end{equation}
Here, of course, the numbers $t(v)$ of allowed absences also need to be adjusted after
the first round. We must reduce $t(v)$ by one, for every $v\in U.$ Using \emph{indicator
functions} $\E_U$ of subsets $U\sb V,$
 \rand\begin{equation}
"\E_U"(v)\,:=\,\begin{cases}
 1&\text{if \,$v\in U,$}\\
 0&\text{if \,$v\in V\sm U,$}
\end{cases}
\end{equation}
we can turn this observation into the following definition\ssp :

\begin{Definition}\label{def.ota}
Let $G=(V,E)$ be a multigraph, $t\DP V\!\to\N$ a function and $m\in\N.$ Being
\emph{on-line \(t\)-avoiding edge \(m\)-colorable} is recursively defined as follows\ssp :
\begin{enumerate}[\((i)\)]
  \item If $m=0,$ then $G$ is on-line \(t\)-avoiding edge \(m\)-colorable if $E=\o.$
  \item If $m>0,$ then $G$ is on-line \(t\)-avoiding edge \(m\)-colorable if for all subsets
      \randd$"U"\sb\supp(t)$ there exists an independent subset \randd$"F"\sb E(G-U)$
      such that $G-F$ is \((t{-}\E_U)\)-avoiding edge \((m{-}1)\)-colorable.
\end{enumerate}
\end{Definition}

To see if a multigraph $G$ is on-line \(t\)-avoiding edge \(m\)-colorable, we have to
apply the recursion step \((ii)\) in our definition up to $m$ times recursively. In each
step, we have to find a suitable response $F$ to every possible choice of $U.$ To keep
track of what happens in these steps, we may denote a $U$ that is chosen in the
$i^{\text{th}}$ step as \randd$"U_i",$ and an $F$ that is chosen in the $i^{\text{th}}$
step as \randd$"F_i"=F_i(U_i).$ Input variables $G=(V,E)$ and $t$ that enter the
$i^{\text{th}}$ round may be denoted as \randd\rand$"G_i"=("V_i,E_i")$ and
\randd$"t_i".$ This notation will make it easier to explain things. For example, it makes
it easy to explain that an on-line scheduling that successfully terminates after $m$
rounds produces an edge \(m\)-coloring $c'\DP E\!\to(m]$ of $G.$ The different $F_i$
just form the color classes of such a \rand$"c'",$ i.e.\
\begin{equation}
\text{$c'(e)=i\ \ \lEqi\ \ e\in F_i.$}
\end{equation}
This, indeed, defines a coloring $E\!\to(m],$ because $F_1\cup F_2\cup\dotsb\cup
F_m=E.$ Otherwise, our recursive scheduling would not have been completed with the
base case $(i)$ of our definition after $m$ rounds. Moreover, the resulting coloring $c'$
avoids the multi-labeling \rand$"c"\DP V\to\cP(\Z^+)$ given by
\begin{equation}
c(v)\,:=\,\{\ssp i\in\Z^+\mit v\in U_i\ssp \},
\end{equation}
because $F_i$ is always chosen outside $U_i\ssp,$ in $G_i-U_i.$ The sets $U_i$
determine the multi-labeling $c$ conversely to how a multi-labeling $c$ has determined
sets $U_i$ in Equation\,\eqref{eq.Ui}\ssp. The only new thing in our recursive definition
is that the $U_i$ are not pre-fixed. Each $U_i$ may be chosen only after the edges
have been partially colored with the colors $1,2,\dotsc,i-1$ already.

We see that our definition actually defines a certain avoiding colorability. Again, the
lowest possible number of colors, i.e.\ the lowest possible number of rounds, is the
crucial parameter for us. We call it the \emph{on-line \(t\)-avoiding chromatic index} of
$G,$ and denoted it by
 \rand\begin{equation}\label{eq.xto}
"\chio^t"(G)\,:=\,\min\{\ssp m\mit\textps{$G$ is on-line \(t\)-avoiding edge
\(m\)-colorable}\ssp \}.
\end{equation}
For $t=0,$ this is again the chromatic index,
\begin{equation}
\chio^0(G)\,=\,\chi^0(G)\,=\,\chi'(G),
\end{equation}
as one easily can check.


\subsection{General Lower Bounds (on-line scheduling)}\label{sec.OLlb}

Time scheduling with complete information is the same as on-line scheduling, but with
pre-fixed sets of absences $U_i.$ Indeed, if a \(t\)-labeling $c\DP V\!\to\cP(\Z^+)$ is
given, we have said that the subset of people that do not appear in round $i$ is
\begin{equation}
U_i\,:=\,\{\ssp v\in V\mit i\in c(v)\ssp \}.
\end{equation}
In our recursive definition this is a valid choice. Indeed, when we enter the $i^\text{th}$
round, we have arrived at $t_i:=t-\E_{U_1}-\E_{U_2}-\dotsb-\E_{U_{i-1}}$ already.
Using that $c$ is a \(t\)-labeling, it is then easy to show that $U_i$ is contained in
$\supp(t_i).$ Therefore, a winning on-line strategy also can be applied to pre-fixed sets
of absences $U_i.$ Hence, an on-line \(t\)-avoiding edge \(m\)-colorable graph $G$ is
\(t\)-avoiding edge \(m\)-colorable. In other words,
\begin{equation}\label{eq.CCO}
\chio^t(G)\,\geq\,\chi^t(G).
\end{equation}
Therefore, all lower bounds for $\chi^t$ remain valid for $\chio^t.$ 
Some of the lower bounds, however, can be improved. As we will see, $\chio^t(G)$ can
be bigger than $\chi^t(G).$ Indeed, if $t$ is constant and $\chi'(G)=m>0,$ then we may
choose
\begin{equation}
U_1=U_2=\dotsb=U_{m-1}\,:=\,\o,
\end{equation}
so that still $t_m=t$ when we enter the $m^\text{th}$ round. At this point also
$E_m\neq\o,$ as $\chi'(G)>m-1.$ Hence, there is a $uv\in E_m.$ The two players $u$
and $v$ may now outmanoeuvre the organizers and never come together for the next
$2t$ rounds. In other words, we may define
\begin{equation}
\begin{split}
U_m=U_{m+1}=\dotsb=U_{m+t-1}\,&:=\,\{u\},\\
U_{m+t}=U_{m+t+1}=\dotsb=U_{m+2t-1}\,&:=\,\{v\}.
\end{split}
\end{equation}
This choice is quite unfortunate for the scheduling, as the edge $uv$ cannot be
selected in $F_m,F_{m+1},\dotsc,F_{m+2t-1}.$ We have to wait until round $m+2t$ to
schedule the match $uv.$ Hence,
\begin{equation}\label{eq.OLlb}
\chio^t(G)\,\geq\,\chi'(G)+2t,
\end{equation}
for all non-edgeless multigraphs $G.$ If $t$ is non-constant, then similarly
\begin{equation}\label{eq.OLlb2}
\chio^t(G)\,\geq\,\chi'(G)+\min_{uv\in E}(t(u)+t(v)).
\end{equation}
For the ordinary $\chi^t(G),$ we had this lower bound just with $\Delta(G)$ in the place
of $\chi'(G).$ For simple graphs $G$ of \emph{class two} (the graphs with
$\chi'(G)>\Delta(G)$) and for constant $t\geq1,$ Lower Bound\,\eqref{eq.OLlb} is
higher than the conjectured value $\Delta(G)+2t$ of $\chi^t(G)$ in
Conjecture\,\ref{con.ub}\ssp. Therefore, this conjecture suggests that, for simple graphs
$G$ and constant $t\geq1,$
\begin{equation}
\chio^t(G)\,>\,\chi^t(G)\quad\text{if}\quad\chi'(G)\,>\,\Delta(G).
\end{equation}

\subsection{General Upper Bounds (on-line scheduling)}\label{sec.OLub}

In this section, we generalize some of the upper bounds for $\chi^t$ to $\chio^t.$
Several of the upper bounds for $\chi^t$ were based on the connection to list edge
colorings. The results in this section can be derived in exactly the same way, but using
the on-line version of list coloring (also called \emph{painting}), a concept that we
introduced in \cite{schPC}\ssp. In this concept, the lists can be modified during an
interactive coloration process. The idea is that, if only positive integers are allowed as
colors, we may use color $1$ at first, of course, only for vertices $v$ whose lists $L_v$
contain color $1.$ Afterwards, before we extend the partial coloring with color $2$,  we
allow changes of the remaining lists $L_v\!\sm\!\{1\}$ that do not change their
cardinalities. This extension process is then repeated with color $3,$ color $4$ and so
forth, where in between, the remaining tails of the color lists may be altered in arbitrary,
possibly unfortunate ways. The related \emph{paintability index} or \emph{on-line list
coloring index} $\chio'$ incorporates the additional flexibility of the color lists. Apart
from that, $\chio'$ is defined exactly as $\chi_\ell'$ in Equation\,\eqref{eq.xl}\ssp. So,
we have $\chio'(G)\leq m$ if and only if $G$ is \emph{edge \(m\)-paintable} (which
means that $G$ can be edge colored from flexible lists of length $m,$ no matter how
these lists are altered after each step of the coloration process). Actually, a recursive
definition of $\chio'\ssp,$ similar to our definition of $\chio^t$ and to Definition\,1.8 in
\cite{schPC}\ssp, would also be possible. Important here is that, even though
$\chio'(G)$ can be bigger than $\chi_\ell'(G),$ the great majority of all list coloring
theorems in graph theory could already be generalized to paintability, see e.g.\
\cite{hks,schPC,schPAT,schPCN}.

Similarly as in the section about $\chi^t,$ our main upper bound is
\begin{equation}
\chio^t(G)\,\leq\,\chio'(G)+2t.
\end{equation}
If the on-line strengthening of the List Edge Coloring Conjecture holds for $G,$ i.e.\ if
$\chio'(G)=\chi'(G),$ then
\begin{equation}
\chio^t(G)\,\leq\,\chi'(G)+2t.
\end{equation}
This would coincide with Lower Bound\,\eqref{eq.OLlb}\ssp. Therefore, as a weakening
of the On-line List Edge Coloring Conjecture, we suspect the following analog to
Conjecture\,\ref{con.ub}\ssp :

\begin{Vermutung}\label{con.ubo}
For every non-edgeless multigraph $G$ and every $t\in\N,$
$$
\chio^t(G)\,=\,\chi'(G)+2t.
$$
\end{Vermutung}

This conjecture is as short as Conjecture\,\ref{con.ub}, but contains the parameter
$\chi'(G)$ which is of more complex nature than $\Delta(G)$. This time, however, we
could allow multiple edges and $t=0$. We also have slightly more evidence, as we will
see later. Moreover, our best general upper bound in Theorem\,\ref{sz.sh} can be
generalized, too. Since Shannon’s Bound and its sharpenings can be generalized to
on-line list edge coloring results, \cite[Theorem\,3.5]{schPC}\ssp, we have the following
improvement of Theorem\,\ref{sz.sh}\ssp :

\begin{Satz}\label{sz.OLsh}
For every multigraph $G=(V,E)$ and every function $t\DP V\!\to\N,$
$$
\chio^t(G)\,\leq\,\max_{uv\in E}\bigl(\max(d(u),d(v))+\cop{\tfrac{1}{2}\min(d(u),d(w))}+t(u)+t(v)\bigr).
$$
In particular, if $t$ is constant,
$$
\chio^t(G)\,\leq\,\cop{\tfrac{3}{2}\Delta(G)}+2t.
$$
\end{Satz}

\subsection{Bipartite Tournaments (on-line scheduling)}\label{sec.OLB}

Since, for bipartite multigraphs $B$
\begin{equation}
\chio'(B)\,=\,\chi_\ell'(B)\,=\,\chi'(B)\,=\,\Delta(B),
\end{equation}
by our strengthening of Galvin's Theorem \cite[Theorem\,3.2]{schPC}\ssp, we can
calculate $\chio^t(B)$ for all $t\in\N.$ We obtain the following sharpening of
Theorem\,\ref{sz.ga}\ssp :

\begin{Satz}\label{sz.gao}
For every non-edgeless bipartite multigraph $B$ and every $t\in\N,$
$$
\chio^t(B)\,=\,\chi^t(B)\,=\,\chi'(B)+2t\,=\,\Delta(B)+2t.
$$
\end{Satz}

As in the section about $\chi^t,$ we can also provide upper bounds for non-constant
$t.$ This is based on our generalization of Borodin, Kostochka and Woodall's Theorem
\cite[Theorem\,3\ssp ]{bkw} to paintability in \cite[Theorem\,3.3\ssp ]{schPC}\ssp.
Exactly as in the section about $\chi^t,$ this yields the following sharpening of
Theorem\,\ref{sz.BKW}\ssp :

\begin{Satz}\label{sz.BKWo}
For every bipartite multigraph $B=(V,E)$ and every function $t\DP V\!\to\N,$
$$
\chio^t(B)\,\leq\,\max_{uv\in E}\bigl(\max(d(u),d(v))+t(u)+t(v)\bigr)\,\leq\,\Delta(B)+\max_{uv\in E}\bigl(t(u)+t(v)\bigr).
$$
\end{Satz}

Together with Lower Bound \eqref{eq.CCO} and Corollary\,\ref{cor.BKW}\ssp, this
implies the following corollary\ssp :

\begin{Korollar}\label{cor.BKWo}
If $B$ is non-edgeless and $t\DP V\!\to\N$ takes the value $t_1$ on one block and the
value $t_2$ on the other block of the bipartition, then
\begin{equation}\label{eq.lbOC2}
\chio^t(B)\,=\,\chi^t(B)\,=\,\chi'(B)+t_1+t_2\,=\,\Delta(B)+t_1+t_2.
\end{equation}
\end{Korollar}

\subsection{Complete Tournaments (on-line scheduling)}\label{sec.OLKn}

In this section, we generalize the results about $\chi^t(K_n).$ 
Since, however, we have the improved Lower Bound\,\eqref{eq.OLlb}\ssp, we have to
be prepared for some differences. H\"{a}ggkvist and Janssen's upper bound still holds
in the on-line version, as we could show in \cite[Theorem\,3.10\ssp ]{schPAT}\ssp, i.e.\
\begin{equation}
\chio'(K_n)\,\leq\,n,
\end{equation}
for all $n\in\Z^+\!.$ This implies that
\begin{equation}
\chio^t(K_n)\,\leq\,n+2t.
\end{equation}
Again, this can be improved by one for all even $n$ for which $\chio'(K_n)=\chi'(K_n),$
and there are partial results in this direction. In \cite{schKp}\ssp, we could prove this
under the additional assumption that $n-1$ is an odd prime. Together with Lower
Bound\,\eqref{eq.OLlb}, this yields the following theorem\ssp :

\begin{Satz}\label{sz.xotKn}
Let $t,n\in\N,$ $n\geq2,$ and assume that $n-1$ is even or prime, then
$$
\chio^t(K_n)\,=\,\begin{cases}
 n+2t-1&\text{if $n$ is even,}\\
 n+2t&\text{if $n$ is odd.}
\end{cases}
$$
\end{Satz}

It is natural to conjecture that this holds for all $n\geq2.$ In fact, this would just be a
small special case of Conjecture\,\ref{con.ubo}\ssp. It also would be just a weakening
of a small special case in the On-line List Edge Coloring Conjecture. That
$\chio^t(K_{10})=10+2t-1$ follows from the fact that the natural action of the symmetric
group $S_{10}$ on the \(1\)-factorizations (edge colorings) of $K_{10}$ has only one
orbit whose length is not divisible by $25,$ \cite[VII.\ssp5.\ssp54\ssp ]{cd}. This fact
shows that, if we count modulo $25,$ the number of \(1\)-factorizations is congruent to
the length of that single orbit. In particular, it is not zero, even if we count the
\(1\)-factorizations in some orbits negative. With that, the ideas in
\cite[Corollary\,3.9\ssp ]{al} or \cite[Theorem\,3.1\ssp ]{schKp} yield
$\chio'(K_{10})=\chi'(K_{10}),$ which implies $\chio^t(K_{10})=10+2t-1.$ Overall,
$n=16$ is the smallest open case in this conjecture (as $n-1=15$ is the first non-prime
after $9,$ and as the case $n=2$ is trivial).


\section{Comparison}\label{sec.comp}

For bipartite multigraphs $B$ the parameters $\Delta(B)$ and $\chi'(B)$ coincide, so
that the new parameters $\chi^t(B)$ and $\chio^t(B)$ also coincide, even if $t$ is just
blockwise constant (by Corollary\,\ref{cor.BKWo}\ssp). For general multigraphs $G,$
however, the situation is different. The parameters $\chi^t(G)$ and $\chio^t(G)$ may
differ, but they seem to differ by at most one if there are no multiple edges. For simple
graphs and if our conjectures hold, there are only two cases (by Vizing's Theorem).
The parameters $\chi^t(G)$ and $\chio^t(G)$ coincide for all \emph{class\,1 graphs}
($\chi'(G)=\Delta(G)$) but differ by one for \emph{class\,2 graphs}
($\chi'(G)=\Delta(G)+1$). Almost all graphs are class\,1 \cite{ew}, but it is usually
difficult to verify that a given graph is actually class\,1 \cite{ho}. Complete graphs
$K_n$ on $n$ vertices are class\,1 if $n$ is even, but class\,2 if $n$ is odd. In
particular, this means that we know $\chi^t(K_n)$ and $\chio^t(K_n)$ for all $n\geq2$
and $t\geq1$, if our conjectures hold. To illustrate what we have learned about
complete tournaments with single absences ($t=1$), we pick up Example\,\ref{ex.K3}
and extend it as follows\ssp:

\begin{Beispiel}\label{ex.Kn}
A group of $n$ chess players wants to arrange a chess tournament. Each player shall
play exactly once against each other player. Moreover, each player shall be allowed to
miss one round. We want to see how many rounds are needed to schedule the full
tournament\ssp:

If the absences do not have to be pre-announced, then $\chio^1(K_n)$ is the best
upper bound on the number of rounds. Based on Theorem\,\ref{sz.xotKn}, the blue
graph in Figure\,\ref{fig.2} shows $\chio^1(K_n)$ as function of $n$ for $n\leq9$.

If the absences are pre-fixed and pre-announced in advance, then $\chio^1(K_n)$ is
the best upper bound on the number of rounds. Based on Theorem\,\ref{sz.x1Kn}, the
red graph in Figure\,\ref{fig.2} shows $\chi^1(K_n)$ as
function of $n$. 
We remind the reader, however, that $\chi^1(K_n)$ is a common upper bound for all
distributions $c\DP(n]\to\Z^+$ of absences. It is best possible only as common upper
bound.

If the absences are actually pre-announced in advance, one may also have a closer
look at their distribution. For complete graphs and $t=1$, we actually could determine
the influence of the distribution $c\DP(n]\to\Z^+$ to the number $\chi'_c(K_n)$ of
rounds that is needed to schedule all games.
Theorem\,\ref{sz.xcKn} tells us precisely how many rounds are needed. 
With the function $c$ as additional parameter, however, it is not possible to print a
meaningful graph for $\chi'_c(K_n)$ in Figure\,\ref{fig.2}. For more theoretical reasons,
it might be interesting to see as brown graph the total chromatic number $\chi_T(K_n)$
(defined in Equation\,\eqref{eq.xT}). It tells the number of rounds if exactly one player is
missing each round. If all players are missing the first round, one gets the chromatic
index plus one $\chi'(K_n)+1$ as number of rounds. We displayed this well studied
graph parameter in green to allow a comparison with our new parameters
$\chi^1(K_n)$ and $\chio^1(K_n)$.\bigskip

\begin{figure}[H]
\begin{center}
\includegraphics[width=.43\textwidth, height=.377\textwidth
]{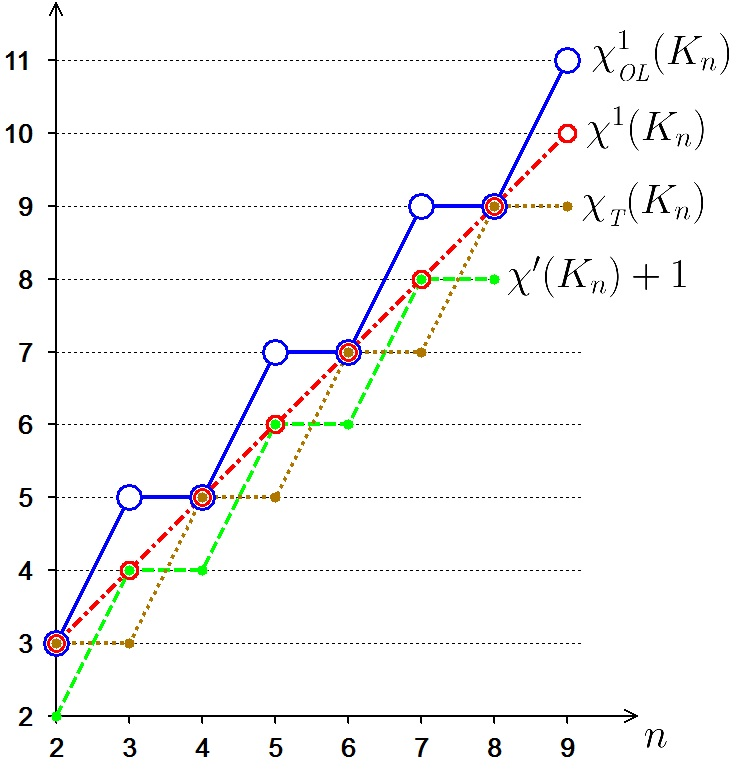}\vspace{-1ex}
\caption{$\chio^1(K_n)$ in comparison with $\chi^1(K_n),$ $\chi_T(K_n)$ and $\chi'(K_n)+1.$}\label{fig.2}
\end{center}\vspace{-3ex}
\end{figure}
\end{Beispiel}

\section{Concluding Remarks}\label{sec.conc}

After the introduction of the parameters $\chi^t(G)$ and $\chio^t(G)$, this paper
focused mainly on Conjecture\,\ref{con.ub} ($\chi^t(G)=\Delta(G)+2t$) and
Conjecture\,\ref{con.ubo} ($\chio^t(G)=\chi'(G)+2t$). Intuitively, these equations look
wrong, as they seem too simple to be true. In particular the first equation appears
unlikely. All commonly studied chromatic numbers are hard to calculate, but, in this
conjecture, the parameter $\chi^t(G)$ can simply be read off from the maximal degree.
This is a serious counterargument. Despite this doubt, however, we could provide
some supporting evidence.

Both conjectures were motivated by our lower bounds, so that the stated values are
certainly not too high. We also were able to prove our conjectures for all bipartite
graphs.
For bipartite graphs, however, both conjectures become the same. So, this particular
case does not provide justification for the finer differences between our conjectures.
We need to consider non-bipartite graphs (and class\,2 graphs) to show that we
adjusted our conjectures in the right way to pre-announced and unannounced
absences.

With respect to unannounced absences and our second conjecture, we showed that
this conjecture holds for all graphs for which the On-Line List Edge Coloring Conjecture
holds. In particular, there are infinitely many complete graphs for which our second
conjecture holds. It would be nice if this result could be extended to all complete
graphs. This, however, would require some new ideas. Eventually,  our approach relies
on \cite[Lemma\,4.1]{schKp}\ssp, and this lemma cannot be generalized in a suitable
way, as unpublished counterexamples show.

With respect to pre-announced absences and our first conjecture, a connection to a
long-standing conjecture was established, too. However, our first conjecture
strengthens this conjecture, the Total Coloring Conjecture, and not the other way
around. Therefore, we cannot simply search for graphs that meet the Total Coloring
Conjecture and then deduce that our first conjecture must also hold for them. Apart
from bipartite graphs and, more generally, all class\,1 graphs for which the On-Line List
Edge Coloring Conjecture holds, we could fully prove our first conjecture only for
complete graphs and $t=1$. We would like to see a proof for all $t\geq1$, preferably
with a full scheduling scheme, like the one inside the proof of
Theorem\,\ref{sz.x1Kn}\ssp. This would be helpful with respect to the organization of
complete tournaments with up to $t$ announced absences per player.

\bigskip
\bigskip

\noindent\textbf{Acknowledgement\ssp:} We thank Jingyi Shi and Zhaoyuan Huang for
the initial calculation of $\chi^1(K_n)$ and $\chio^1(K_n)$ for some small $n.$



\begin{thebibliography}{10}
\bibitem[Al]{al}
  N.\,Alon: \textit{Restricted Colorings of Graphs.}\\
  In: Surveys in combinatorics, 1993.
  London Math.\ Soc.\ Lecture Notes Ser.\ 187,\\
  Cambridge Univ.\ Press, Cambridge 1993, 1–33.
\bibitem[An]{an}
  L.\,D.\,Andersen: \textit{Embedding Latin Squares with Prescribed Diagonal.}\\
  Annals of  Discrete Mathematics 15 (1982), 9–26.
\bibitem[BKW]{bkw}
  O.\,V.\,Borodin, A.\,V.\,Kostochka, D.\,R.\,Woodall:\\
  \emph{List Edge and List Total Colourings of Multigraphs.}\\
  J.\ Combin.\ Theory Ser.\ B 71(2) (1997), 184–204.
\bibitem[CoDi]{cd}
  C.\,J.\,Colbourn, J.\,H.\,Dinitz: \textit{The CRC Handbook of Combinatorial Designs.}\\
  Sec.\ Edition (Discr.\ Mathematics and Its Applications)
  Chapman and Hall/CRC 2006.
\bibitem[Cr]{cr}
  A.\,B.\,Cruse: \textit{On Embedding Incomplete Symmetric Latin Squares.}\\
  J.\ Combin.\ Theory Ser.\ A 16 (1974), 18–27.
\bibitem[Di]{di}
  R.\,Diestel: \textit{Graph Theory (4rd edition).}
  Springer, Berlin 2010.
\bibitem[DrKn]{dk}
  A.\,Drexl and S.\,Knust:\\ \textit{Sports league scheduling: Graph-\ and Resource-Based Models.}\\
  Omega, vol.\,35, issue\,5 (2007), 465–471.
\bibitem[ErWi]{ew} P.~Erd\H{o}s and R.~J. Wilson.
  \textit{On the chromatic index of almost all graphs.}\\
  Journal of Combinatorial Theory, Series B, 23(2-3):255--257, 1977.
\bibitem[FiWi]{fiwi}
    S.\,Fiorini and R.\,J.\,Wilson: \textit{Edge-Colourings of Graphs.}\\
    Research Notes in Mathematics, Pitman, 1977.
\bibitem[Ga]{ga}
    F.\,Galvin: \textit{The List Chromatic Index of a Bipartite Multigraph.}\\
    J.\,Combin.\ Theory Ser.\ B 63 (1995), 153–158.
\bibitem[H\"{a}Ja]{haja}
    R.\,H\"{a}ggkvist  and J.\,Janssen: \textit{New Bounds on the List-Chromatic Index
    of the Complete Graph and Other Simple Graphs.}\\
    Combin.\ Probab.\ Comput.\ 6 (1997), 295–313.
\bibitem[HKS]{hks}
  J.\,Hladk\'{y}, D.\,Kr\'{a}l', U.\,Schauz: \textit{Brooks' Theorem via the Alon-Tarsi Theorem.}\\
  Discrete Mathematics 310 (2010), 3426–3428.
\bibitem[Ho]{ho}
  I.\,Holyer: \textit{The NP-completeness of edge-coloring.}\\
  SIAM Journal on Computing, 10(4):718--720, 1981
\bibitem[JeTo]{jeto}
  T.\,R.\,Jensen, B.\,Toft: \textit{Graph Coloring Problems.} Wiley, New York 1995.
\bibitem[KKRU]{kkru}
  G.\,Kendall, S.\,Knust, C.\,Ribeiro, S.\,Urrutia: \textit{Scheduling in Sports:\\ An Annotated Bibliography.}\
    Computers \& Operations Research 37 (2010), 1–19.
\bibitem[Kn]{kn}
  S.\,Knust: \textit{Classification of Literature on Sports Scheduling.}\\
  http://www2.informatik.uni-osnabrueck.de/knust/sportssched/sportlit\_class/\\ (assesed 1.May 2016)
\bibitem[LeTh]{lt}
  R.\,Lewis, J.\,Thompson: \textit{On the Application of Graph Colouring
    Techniques in Round-Robin Sports Scheduling.}\\
    Computers \& Operations Research 38 (2011), 190–204.
\bibitem[RaTr]{rt}
  R.\,V.\,Rasmussen, M.\,A.\,Trick:
  \textit{Round Robin Scheduling – A Survey.}\\
  European Journal of Operational Research 188 (2008), 617–636.
\bibitem[Sch2]{schPC}
  U.\,Schauz: \textit{Mr.\ Paint and Mrs.\ Correct.}\\
  The Electronic Journal of Combinatorics 15 (2008), \#R145.
\bibitem[Sch3]{schPAT}
  U.\,Schauz: \textit{Flexible Color Lists in Alon and Tarsi's Theorem,\\
  and Time Scheduling with Unreliable Participants.}\\
  The Electronic Journal of Combinatorics 17/1 (2010), \#R13.
\bibitem[Sch4]{schPCN}
  U.\,Schauz: \textit{A Paintability Version of the Combinatorial Nullstellensatz,\\
  and List Colorings of \(k\)-partite \(k\)-uniform Hypergraphs.}\\
  The Electronic Journal of Combinatorics 17/1 (2010), \#R176.
\bibitem[Sch5]{schKp}
  U.\,Schauz:\\ \textit{Proof of the List Edge Coloring Conjecture for Complete Graphs of Prime Degree.}\\
  The Electronic Journal of Combinatorics 21/3 (2014), \#P3.43.
\bibitem[Ya]{ya} H.\,P.\,Yap: \emph{Total Colourings of Graphs.}\\
  Lect.\ Notes in Math.\ 1623, Springer, Berlin 1996.
\end{thebibliography}

\end{document}